\documentclass[11pt]{article}
\input amssym.def
\input amssym

\begin{document}
\bibliographystyle{alpha}

\newtheorem{thm}{Theorem}
\newtheorem{defin}[thm]{Definition}
\newtheorem{lemma}[thm]{Lemma}
\newtheorem{propo}[thm]{Proposition}
\newtheorem{cor}[thm]{Corollary}
\newtheorem{conj}[thm]{Conjecture}
\newtheorem{exa}[thm]{Example}

\begin{center}
{\LARGE  \bf Measurable Categories and 2-Groups}
\vspace*{1cm}

\centerline{\parbox{2.5in}{Louis Crane and David N. Yetter\\
Department of Mathematics\\
Kansas State University\\
Manhattan, KS  66506 }}

\end{center}
\vspace*{1cm}

\noindent{\bf Abstract:}  Using the theory of measurable categories developped in
\cite{Y.meas}, we provide a notion of representations of 2-groups more
well-suited to physically and geometrically interesting examples than that
using 2-VECT (cf. \cite{KV}).  Using this theory we sketch a 2-categorical
approach to the state-sum model for Lorentzian quantum gravity proposed in
\cite{CY.Lorentz}, and suggest state-integral constructions for 4-manifold 
invariants.

\section{Introduction}

In Section \ref{meascat} we recall the relevant definitions and results from
\cite{Y.meas}.  In Section \ref{2groups} we recall the relevant definitions
and well-known results concerning 2-groups, discuss examples of particular
importance for physics and group representation theory including the 
Poincar\'{e}
2-group first considered by Baez \cite{B}.  Section \ref{repth} develops the
representation theory of 2-groups in the 2-category ${\bf Meas}$.
Section\ref{ss} outlines the common features of anticipated
constructions both of a large family of new topological quantum 
field theories arising from the general theory and of a 
new model for quantum general relativity arising from the representation
theory of the Poincar\'{e} 2-group.

\section{Measurable categories} \label{meascat}

Throughout this work, we assume that all Borel spaces have measurable points.  
This restriction is necessary for the definition of ``measurable functors'' 
given below.  However, it excludes only pathological examples:  Euclidean 
space with either Borel or Lebesgue measure, the Borel structure on any 
locally compact group, and all discrete Borel spaces satisfy the condition.

We begin by recalling several definitions from \cite{Y.meas}

\begin{defin}
A {\em measurable field of Hilbert spaces} ${\cal H}$ on a Borel space
$(X,S)$ is a pair $({\cal H}_x, {\cal M}_{\cal H})$, where ${\cal H}_x$ is
an X-indexed family of Hilbert spaces, and ${\cal M}_{\cal H} = {\cal M}$ 
is a linear
subspace of $\prod_{x \in X} {\cal H}_x$ (the product as vector-spaces)
satisfying

\begin{enumerate}
\item $\forall \xi \in {\cal M} \;  x \mapsto \|\xi(x)\|_x$ is measurable
\item $\forall \eta \in \prod_{x \in X} {\cal H}_x \; 
	x \mapsto \langle \eta(x) | \xi(x) \rangle_x$ is measurable for
	all $\xi \in {\cal M}$ implies $\eta \in {\cal M}$
\item $\exists \{ \xi_i \}_{i = 1}^\infty \subset {\cal M}$ such that
	${\xi_i(x)}_{i = 1}^\infty$ is dense in ${\cal H}_x$ for all $x \in X$
\end{enumerate}

An {\em almost measurable field of Hilbert spaces} ${\cal H}$ on a Borel
space $(X,S)$ is a pair $({\cal H}_x, {\cal M}_{\cal H})$ as above, 
satisfying 1. and 2., but not necessarily 3.

\end{defin}

\begin{defin}
A {\em measurable field of bounded operators} $\phi$ from ${\cal H}$ to 
${\cal K}$, for ${\cal H}$ and ${\cal K}$ (almost) measurable fields of 
Hilbert spaces is an
X-indexed family of bounded operators $\phi_x \in B({\cal H}_x, {\cal K}_x)$
such that $\xi \in {\cal M}_{\cal H}$ implies 
$\phi(\xi) \in {\cal M}_{\cal K}$,
where $\phi(\xi)_x = \phi_x(\xi_x)$.

A measurable field of bounded operators is {\em bounded} if
$x \mapsto \|\phi_x\|_x$ is a bounded real-valued function  (Here $\| \; \|_x$ denotes the
operator norm on $B({\cal H}_x, {\cal K}_x)$.)
\end{defin}

We can then organize these into a category:

\begin{defin}
The {\em category of measurable fields of Hilbert spaces on} $(X,S)$ has as
objects all measurable fields of Hilbert spaces on $(X,S)$ and as arrows
all bounded fields of bounded operators on $X$.  Source, target, 
identity arrow and composition are obvious.  We denote this category by
${\bf Meas}(X,S)$.

Similarly, the {\em category of almost measurable fields of Hilbert spaces} on
$(X,S)$ has as objects all almost measurable fields of Hilbert spaces on 
$(X,S)$ and as arrows all bounded fields of bounded operator 
between them.  We denote this category by ${\bf AlMeas}(X,S)$.
\end{defin}

As in \cite{Y.meas} these categories in turn are organized into a 
2-category.  

In describing the representations of 2-groups we need only
invertible functors, we do not at first need
to consider the full theory of measurable
functors developed in \cite{Y.meas}. 
In describing the representations it is a matter of indifference whether we
work in the 2-category {\bf Meas} of \cite{Y.meas}, with
all measurable functors as 1-arrows, or the
2-category ${\bf Meas}_{add}$, with ${\Bbb C}$-linear additive
functors as arrows.

Recall from \cite{Y.meas}

\begin{thm}
Any additive functor with an additive inverse $\Phi:{\bf
Meas}(X,S)\rightarrow {\bf Meas}(Y,T)$ is equivalent to a
functor induced by pullback along an
invertible measurable transformation $\tilde{\Phi}:Y\rightarrow X$.
\end{thm}

\noindent and

\begin{thm}
Any ${\Bbb C}$-linear natural endomorphism of an invertible
additive functor is given by fiberwise scalar multiplication by
an bounded scalar valued function.  2-dimensional composition is 
given by multiplication of functions, and 1-dimensional composition with 
1-arrows is given by pullback.
\end{thm}

In describing the intertwiners, however, it will be necessary to choose
either {\bf Meas} or ${\bf Meas}_{add}$.  We prefer the former, as its
structure is better understood.

For any measurable field of Hilbert spaces $\cal K$ on $X \times Y$
and a $Y$-indexed family of measures on $X$, $\{ \mu_y\}$, let

\[ \Phi_{{\cal K},\{\mu_y \}}({\cal H})_y =
	\int^\oplus_X {\cal H}_x \otimes {\cal K}_{<x,y>} d\mu_y(x) \]

\noindent with ${\cal M}_{\Phi({\cal H})}$ given as the closure
under condition 2 of the set $\{ \int^\oplus_{A_y} \eta(x) \otimes \kappa(x)
d\mu_y(x) | \eta \in {\cal M}_{\cal H}; \kappa \in {\cal M}_{\cal K} \}$.

\begin{defin}
A functor from ${\bf AlMeas}(X,S)$ to ${\bf AlMeas}(Y,T)$ is
{\em measurable} if it ${\Bbb C}$-linear
equivalent to one of the following form
$\Phi_{{\cal K},\{\mu_y \}}$

A {\em measurable functor} from ${\bf Meas}(X,S)$ to ${\bf AlMeas}(Y,T)$
is the restriction of a measurable functor from ${\bf AlMeas}(X,S)$
to ${\bf AlMeas}(Y,T)$, while a {\em measurable functor} from ${\bf Meas}(X,S)$
to ${\bf Meas}(Y,T)$ is the factorization of a measurable functor from
${\bf Meas}(X,S)$
to ${\bf AlMeas}(Y,T)$ through ${\bf Meas}(Y,T)$, provided it admits
such a factorization.
\end{defin}

\section{2-Groups} \label{2groups}

\begin{defin}
A {\em categorical group} is a group object in the category of (small) 
categories.
\end{defin}

It turns out that a categorical group is necessarily a groupoid:  an amusing 
little
exercise using the middle-four interchange law for the functoriality of the 
group
law and the covariant(!) functor $^{-1}$ shows that every arrow has an inverse.

A moment's consideration reveals that categorical groups are in fact strict 
monoidal
categories equipped with a very strong type of two-sided dual.  Because of 
this, we denote the identity object by $I$ as is customary in monoidal
categories. We may thus
use the usual trick to regard categorical groups as 2-categories with a 
single object
(and their 'objects' as 1-arrows, 'arrows' as 2-arrows, the group law as
1-dimensional composition, and the composition as 2-dimensional composition).
When we do this, we refer to the resulting 2-category as a {\em 2-group}.

It is also easy to see that the group law on any categorical group $\frak C$
induces group structures on $Ob({\frak C})$ and $Arr({\frak C})$.

Categorical groups have been studied previously as models of homotopy 2-types
(cf. \cite{BS} \cite{Y.homotopy2}). Brown and Spencer 
\cite{BS} show
that a categorical group is equivalent to a {\em crossed module}:

\begin{defin}
A {\em crossed module} is a homomorphism of groups $\partial :E\rightarrow G$
together with an action $\triangleright$ of $G$ on $E$ by automorphisms, such that

\[ \partial (g\triangleright e) = g(\partial e) g^{-1} \]

\[ (\partial e)\triangleright \varepsilon = e\varepsilon e^{-1} \]
\end{defin}

The group $G$ is called the {\em base group}, while the group $E$ is called 
the {\em fiber group} (other authors call $E$ ``the principal group'', 
but we prefer the name fiber group to emphasise a similarity with fiber 
bundles).  

The equivalence given by Brown and Spencer \cite{BS} 
arises as follows:  Given a 
categorical group $\frak G$ we let $G = Ob({\frak G})$, the group of arrows 
of $\frak G$, and $E \subset Arr({\frak G})$ be the group of all arrow with 
$I$, the identity object, as source.  $\partial$ is then the restriction of 
the target map to $E$, while the action of $G$ on $E$ is given by conjugation 
(in the group of arrows) by the (identity arrow of) the object.

Conversely, given a crossed module, one can form a categorical group by taking
$G$ as the group of objects, and the semidirect product of $E$ and $G$ with 
product $(e,g)(\varepsilon, \gamma) = (eg\triangleright \varepsilon, g\gamma)$
as group
of arrows.  Source and target are given by ${\rm source}(e,g) = g$, 
${\rm target}(e,g) = g\partial(e)$; the identity arrow on $g \in G$ is
$(1_I,g)$, and composition is given by $(e,g)(f,h) = (e\dot f, g)$ whenever
$h = g\partial(e)$.  

The reader is left to complete the proof or refered to \cite{BS}.

We will be particularly interested in categorical groups with a somewhat 
simpler structure:

\begin{defin} 
A categorical group is {\em automorphic} if all of its arrows are 
automorphisms.
\end{defin}

In the crossed module picture, this is equivalent to the map $\partial$ being 
the trivial homomorphism $I$.  Also observe that in this case the fiber group 
is necessarily abelian:  the categorical composition provides a group law for 
which the group law in the crossed module structure is given by homomorphism.
It then follows by the theorem of Eckmann-Hilton \cite{EH} that the two group
laws coincide and are abelian.  We will apply the adjective ``automorphic'' 
to the structure whether we consider it as a categorical group, a crossed 
module, or a 2-group.

We can use group representations to construct automorphic categorical groups, 
including the example closely related to conjectural constructions for 
lorentzian quantum gravity, the Poincar\'{e} 2-group.

\begin{defin}
Let $G$ be a group, and $V, \rho$ be a representation of $G$, the 
{\em representational categorical group} $\rho \natural G$ is the automorphic 
categorical group given as a crossed module with base group $G$ and fiber
group $V,+$ by the action $g\triangleright v = \rho(g)(v)$.

In the case of $G = SO(3,1)$, with ${\Bbb R}^4, \rho$ the natural
action of $SO(3,1)$ by rotations (and boosts) on Minkowski space, the
representational categorical group is the Poincar\'{e} 2-group of \cite{B},
which we denote $\frak P$.
\end{defin}

\section{The Representation Theory of 2-groups} \label{repth}

As with groups, the study of representations means the study of maps to 
particularly nice or well understood examples.  Early attempts
(cf. Barret and Mackaay \cite{BM}) to represent 
2-groups in one or another of the versions of the category $2-VECT$ first 
describe by Kapranov and Voevodsky \cite{KV} have foundered on the paucity of 
examples:  with the exception of categorical groups with (pro)finite groups 
of objects, there will be too few¹ representations to provide a 
satisfying theory.   Indeed, representations in $2-VECT$ cannot be 
collectively faithful unless the group of objects is profinite.

The theory of measurable categories was developed in \cite{Y.meas}, 
motivated by the weak analogy to similar problems in the representation 
theory of non-compact groups, precisely to overcome this difficulty.

\begin{defin}
A {\em (measurable) representation} of a 2-group $\frak G$ is a 2-functor
$R:{\frak G}\rightarrow {\bf Meas}$.  A {\em 1-intertwiner} between
two representations $R$ and $R^\prime$ is a 2-natural transformation 
$\phi:R\rightarrow R^\prime$.  A {\em 2-intertwiner} between two (parallel) 
1-intertwiners $\phi:R\rightarrow R^\prime$ and $\psi:R\rightarrow R^\prime$ 
is a modification $m:\phi:\rightarrow \psi$.
\end{defin}

As with representations of groups, the task of understanding the 
representations of a 2-group consists primarily in decomposing arbitrary 
representations into simpler ones, and understanding the irreducible or
indecomposable examples.  Before this can be attempted, however, it is 
necessary to unwind the previous definition to identify the structures 
involved in terms of more familiar group theoretic and representation 
theoretic notions.

In the present work, we confine ourselves to the consideration of automorphic 
2-groups.  

Now for any 2-group $\cal G$, a representation $R$ assigns to the unique 
object a category $R(*) = {\bf Meas}(X_R,S_R)$ for some measureable space
$(X_R,S_R)$, and to each 1-arrow $G$ an invertible measurable functor
$R(G):{\bf Meas}(X_R,S_R)\rightarrow {\bf Meas}(X_R,S_R)$.  Now, by a result 
of Yetter \cite{Y.meas} any such functor is induced by an invertible 
measurable transformation $R_G:(X_R,S_R)\rightarrow (X_R,S_R)$.  Thus at the 
level of 1-arrows, switching to the language of crossed modules, we may say 
that a representation of a 2-group is specified by a measurable action of 
the base group on a Borel space.

In the case of an automorphic 2-group, any 2-arrow has the same source and
target, and is thus of the form $g:G\rightarrow G$.  The image is thus a
natural automorphism of the functor $R(G)$.  1-dimensional composition with
$G^{-1}$ reduces the problem of describing these to describing the natural
automorphisms of the identity functor on ${\bf Meas}(X_R,S_R)$.  Now by 
the result of \cite{Y.meas} any natural transformation from the identity
functor to itself is determined by an bounded scalar valued 
function.  The preservation of 2-dimensional composition imposes the 
condition that, 
for $\eta:I\rightarrow I$, the
values $R(\eta)(x)$ 
for each $x \in X$ form a character of the fiber group $E$, 
while preservation of the 1-dimensional composition of 2-arrows with 1-arrows 
imposes the condition that $R(F\circ g):R(F\circ G)\rightarrow R(F\circ G)$ 
(resp. $R(g\circ H):R(G\circ H)\rightarrow R(G\circ H)$) is the left 
(resp. right) translation of the function $R(g)(x)$ by $R(F)$ (resp. $R(H)$).

To understand this better, recall that the 2-arrows of an automorphic
2-group form a group under the 1-dimensional composition $\circ$ which
is isomorphic to the semidirect product of $E$ and $G$.  
We may thus index them by pairs
$(\eta, G)$, where $\eta:I\rightarrow I$.  In these terms, the 1- and
2- dimensional compositions are given by $(\eta, G)\circ (\epsilon, H)
= (\eta G\triangleright \epsilon, GH)$ and $(\eta, G)(\eta^\prime, G) 
= (\eta \eta^\prime, G)$ respectively. 

It is an amusing exercise to verify directly in terms of this formulation 
that the middle-four interchange
law holds.

 In terms of this indexing, what
we denoted $R(\eta)(x)$ above is $R(\eta, I)(x)$.

The images of all the 2-arrows may be identified with bounded functions
$R(\eta, G)(x)$ which multiply the fiber at $x$ before the action of
$R(G)$ displaces it.  With this identification, the preservation of
the 1-dimensional composition becomes

\[ R(\eta ,G)(x) R(\epsilon, H)(G(x)) = R(\eta G\triangleright \epsilon, GH)(x)
\]

\noindent  Observe that for any 1-arrow $G$, and $e$ the identity
in the fiber group, we have $R(e,G)(x) \equiv 1$.  This fact, together
with the condition above, gives

\[ R(\eta, G)(x) = R(\eta, G)(x)R(e, H)(G(x)) = R(\eta, GH)(x) \]

\noindent and

\[ R(\epsilon, H)(G(x)) = R(e, G)(x)R(\epsilon, H)(G(x)) = R(G\triangleright
\epsilon, GH)(x) \]

From the first, we see that the bounded function representing 2-arrows
$(\eta, G)$ is independent of $G$.  Thus, we may let ${\cal R}(\eta)(x) =
R(\eta, G)(x)$ for any (and thus all) $G$.

In terms of this, the second condition becomes 

\[ {\cal R}(\epsilon)(G(x)) = {\cal R}(G\triangleright \epsilon)(x) \]

A simple calculation shows that this condition, together with the 
preservation of 2-dimensional composition implies the preservation
of 1-dimensional composition.

We have thus obtained the following:

\begin{propo}
A representation of an automorphic 2-group with base group $\cal G$ and
fiber group $E$ is given by a measurable action of $\cal G$ on a 
measurable space $(X,S)$, together with an $X$-indexed family 
${\cal R}(\epsilon)(x)$ of 
characters of $E$, which is $\cal G$-equivariant in the sense that

\[ {\cal R}(\epsilon)(G(x)) = {\cal R}(G\triangleright \epsilon)(x) \]

\end{propo}

Note that if the action of the base group $G$ on the Borel space $X$ is transitive,
then choosing a point of $x \in X$ gives an indexing of the representations
by characters of the fiber group assigned to $x$.  If the action base
group
is also free, the representations with the given action of the base group
are indexing is by all the characters of the fiber group.  Otherwise, the
indexing is by the characters fixed by the stablizer of $x$.

This indexing, is, however, dependent on the necessarily non-canonical
choice of a point in the $G$-equivariant Borel space $X$.

In the case of a non-transitive action, choosing a point in each orbit gives
a indexing of the representations with the given action of the base group by
maps from the orbits to the set of characters of the fiber group.

We now turn to the task of giving a similar description to the 1-intertwiners
for representations of an automorphic 2-group.

Fix two representations $R$ and $S$, and let $X$ and $Y$ respectively be the
measure space to which each maps the unique object.

Now, a 1-intertwiner $\Phi$ is a 2-natural transformation.  It thus assigns to
the unique object $*$ a measurable functor $\Phi_*$.  As we will only require
a description up to natural equivalence, we may assume that
$\Phi_* = \Phi_{{\cal K}, \{ d\mu_y \}}$ for some $\cal K$, a measurable field of Hilbert spaces on $X\times Y$, and
$d\mu_y$ measures on $X$.

To each 1-arrow $G$ (element of the base group in the crossed module picture), it
assigns a filler for the square

\begin{center}

\setlength{\unitlength}{4144sp}%
\begingroup\makeatletter\ifx\SetFigFont\undefined
% extract first six characters in \fmtname
\def\x#1#2#3#4#5#6#7\relax{\def\x{#1#2#3#4#5#6}}%
\expandafter\x\fmtname xxxxxx\relax \def\y{splain}%
\ifx\x\y   % LaTeX or SliTeX?
\gdef\SetFigFont#1#2#3{%
  \ifnum #1<17\tiny\else \ifnum #1<20\small\else
  \ifnum #1<24\normalsize\else \ifnum #1<29\large\else
  \ifnum #1<34\Large\else \ifnum #1<41\LARGE\else
     \huge\fi\fi\fi\fi\fi\fi
  \csname #3\endcsname}%
\else
\gdef\SetFigFont#1#2#3{\begingroup
  \count@#1\relax \ifnum 25<\count@\count@25\fi
  \def\x{\endgroup\@setsize\SetFigFont{#2pt}}%
  \expandafter\x
    \csname \romannumeral\the\count@ pt\expandafter\endcsname
    \csname @\romannumeral\the\count@ pt\endcsname
  \csname #3\endcsname}%
\fi
\fi\endgroup
\begin{picture}(2880,2385)(46,-1636)
\thinlines
\put(496,344){\vector( 0,-1){1440}}

\put(631,479){\vector( 1, 0){1980}}

\put(2746,344){\vector( 0,-1){1440}}

\put(586,-1321){\vector( 1, 0){1980}}

\put(451,389){\makebox(0,0)[lb]{\smash{\SetFigFont{12}{14.4}{it}$X$
}}}
\put(2701,389){\makebox(0,0)[lb]{\smash{\SetFigFont{12}{14.4}{it}$Y$
}}}
\put(451,-1411){\makebox(0,0)[lb]{\smash{\SetFigFont{12}{14.4}{it}$X$
}}}
\put(2701,-1411){\makebox(0,0)[lb]{\smash{\SetFigFont{12}{14.4}{it}$Y$
}}}
\put( 46,-421){\makebox(0,0)[lb]{\smash{\SetFigFont{12}{14.4}{it}$R(G)$
}}}
\put(2926,-421){\makebox(0,0)[lb]{\smash{\SetFigFont{12}{14.4}{it}$S(G)$
}}}
\put(1351,614){\makebox(0,0)[lb]{\smash{\SetFigFont{12}{14.4}{it}$\Phi_*$
}}}
\put(1441,-1636){\makebox(0,0)[lb]{\smash{\SetFigFont{12}{14.4}{it}$\Phi_*$
}}}
\put(1509,-518){\makebox(0,0)[lb]{\smash{\SetFigFont{12}{14.4}{it}$\Phi_G$
}}}
\put(1433,-323){\makebox(0,0)[lb]{\smash{\SetFigFont{12}{14.4}{it}$\Rightarrow$
}}}
\end{picture}

\end{center}

\noindent so that the composition of 1-arrows (the group law in the base group)
is carried to pasting composition, and for each 2-arrow $\eta:G\rightarrow G$,
the "pillow"

\begin{center}

\setlength{\unitlength}{4144sp}%
\begingroup\makeatletter\ifx\SetFigFont\undefined
% extract first six characters in \fmtname
\def\x#1#2#3#4#5#6#7\relax{\def\x{#1#2#3#4#5#6}}%
\expandafter\x\fmtname xxxxxx\relax \def\y{splain}%
\ifx\x\y   % LaTeX or SliTeX?
\gdef\SetFigFont#1#2#3{%
  \ifnum #1<17\tiny\else \ifnum #1<20\small\else
  \ifnum #1<24\normalsize\else \ifnum #1<29\large\else
  \ifnum #1<34\Large\else \ifnum #1<41\LARGE\else
     \huge\fi\fi\fi\fi\fi\fi
  \csname #3\endcsname}%
\else
\gdef\SetFigFont#1#2#3{\begingroup
  \count@#1\relax \ifnum 25<\count@\count@25\fi
  \def\x{\endgroup\@setsize\SetFigFont{#2pt}}%
  \expandafter\x
    \csname \romannumeral\the\count@ pt\expandafter\endcsname
    \csname @\romannumeral\the\count@ pt\endcsname
  \csname #3\endcsname}%
\fi
\fi\endgroup
\begin{picture}(2639,1905)(144,-1418)
\thinlines
\put(608,-474){\oval(360,1216)[tl]}
\put(608,-474){\oval(360,1214)[bl]}
\put(608,-1081){\vector( 1, 0){0}}

\put(886,-1081){\vector(-1, 0){0}}
\put(886,-474){\oval(360,1214)[br]}
\put(886,-474){\oval(360,1216)[tr]}

\put(2251,-459){\oval(360,1216)[tl]}
\put(2251,-459){\oval(360,1214)[bl]}
\put(2251,-1066){\vector( 1, 0){0}}

\put(2513,-1066){\vector(-1, 0){0}}
\put(2513,-459){\oval(360,1214)[br]}
\put(2513,-459){\oval(360,1216)[tr]}

\put(954,217){\vector( 1, 0){1215}}

\put(901,-1148){\vector( 1, 0){1215}}

\put(2783,-826){\makebox(0,0)[lb]{\smash{\SetFigFont{12}{14.4}{it}$S(G)
$}}}
\put(586,-271){\makebox(0,0)[lb]{\smash{\SetFigFont{12}{14.4}{it}$R(\eta)
$}}}
\put(1111,-938){\makebox(0,0)[lb]{\smash{\SetFigFont{12}{14.4}{it}$R(G)
$}}}
\put(684,119){\makebox(0,0)[lb]{\smash{\SetFigFont{12}{14.4}{it}$X
$}}}
\put(684,-1208){\makebox(0,0)[lb]{\smash{\SetFigFont{12}{14.4}{it}$X
$}}}
\put(1455,-1418){\makebox(0,0)[lb]{\smash{\SetFigFont{12}{14.4}{it}$\Phi_*
$}}}
\put(1440,352){\makebox(0,0)[lb]{\smash{\SetFigFont{12}{14.4}{it}$\Phi_*
$}}}
\put(1808,-15){\makebox(0,0)[lb]{\smash{\SetFigFont{12}{14.4}{it}$S(G)
$}}}
\put(144,-39){\makebox(0,0)[lb]{\smash{\SetFigFont{12}{14.4}{it}$R(G)
$}}}
\put(2214,-586){\makebox(0,0)[lb]{\smash{\SetFigFont{12}{14.4}{it}$S(\eta)
$}}}
\put(2319,-1216){\makebox(0,0)[lb]{\smash{\SetFigFont{12}{14.4}{it}$Y
$}}}
\put(1269,-338){\makebox(0,0)[lb]{\smash{\SetFigFont{12}{14.4}{it}$\Phi_G
$}}}
\put(1718,-826){\makebox(0,0)[lb]{\smash{\SetFigFont{12}{14.4}{it}$\Phi_G
$}}}
\put(2311,157){\makebox(0,0)[lb]{\smash{\SetFigFont{12}{14.4}{it}$Y
$}}}
\put(623,-473){\makebox(0,0)[lb]{\smash{\SetFigFont{12}{14.4}{it}$\Rightarrow
$}}}
\put(1201,-158){\makebox(0,0)[lb]{\smash{\SetFigFont{12}{14.4}{it}$\Rightarrow
$}}}
\put(1651,-638){\makebox(0,0)[lb]{\smash{\SetFigFont{12}{14.4}{it}$\Rightarrow
$}}}
\put(2258,-391){\makebox(0,0)[lb]{\smash{\SetFigFont{12}{14.4}{it}$\Rightarrow
$}}}
\end{picture}

\end{center}

\noindent commutes.

As with the representations themselves, we must now unwind this definition to
discover what more it entails in more familiar terms.  Using the additivity
properties shown in \cite{Y.meas}, it suffices to consider the components
of the filler $\Phi$ at 
partial measurable line bundles.  The action of the base group (1-arrows) provides
a restriction on the partial measurable line bundles:  they must be equivariant
under the action, and thus supported on orbits.

The filler $\Phi_G$ thus has as components bounded
fields of operators

\[ \Phi_{G,y}({\cal H}): \int^\oplus_{x\in X} {\cal H}_{R(G)(x)}\otimes
{\cal K}_{<x,y>} d\mu_y(x)
\rightarrow \int^\oplus_{x \in X} {\cal H}_x \otimes {\cal K}_{<x,S(G)(y)>} 
d\mu_{S(G)(y)}(x) \]

The passage of the composition of 1-arrows to pasting composition is then
given by the condition

\[ \Phi_{GH,y} = \Phi_{G,y}\Phi_{H,S(G)(y)} \]

The commutativity of the pillow for $\eta:G\rightarrow G$
becomes

\[ \Phi_{G,y} \cdot {\cal S}(\eta)(S(G)(y)) = \int^\oplus \cdot{\cal R}(\eta)(R(G)(x)) 
d\mu_y(x) \Phi_{G,y} \]

\noindent but, since multiplication by a field of scalars commutes with any field
of bounded operators, this becomes the condition that $\Phi_{G,y}$ coequalizes
$\cdot {\cal S}(\eta)(S(G)(y))$ and $\int^\oplus \cdot{\cal R}(\eta)(R(G)(x) 
d\mu_y(x)$, and thus, since the $\Phi_{G,y}$ is invertible that they are equal.

It therefore follows for all $\eta:G\rightarrow G$ and all $y\in Y$ 
that ${\cal R}(\eta)(R(G)(x)) ={\cal S}(\eta)(S(G)(y))$
$\mu_y$-almost-everywhere.

However, the condition on the fields of characters which define the 2-arrow
part of a representation allows us to replace this condition with
${\cal R}(G\triangleright \eta)(x) = {\cal S}(G\triangleright \eta)(x)$
$\mu_y$-a.e.
But, since this must hold for all $G$ and $\eta$, if we replace
$\eta$ with $G^{-1}\triangleright \eta$, the condition reduces to
${\cal R}(\eta)(x) = {\cal S}(\eta)(y)$  $\mu_y$-a.e.

Consideration of the direct integral of the source, target, and fields
of operators $\Phi_{G,y}$ with respect to any measure on $Y$ equivariant
with respect to the action of the 1-arrows given by $S$ shows that any
such direct integral carries a representation of the group of 2-arrows
under 1-dimensional composition.  But, it is a representation with 
additional structure:  The base group acts by operators which are
the composition of a diagonal operator with the translation operator
on the underlying measure space defined by $S$.  The fiber group acts
by multiplication by the equivariant field of characters given by $S$.
And, finally, the spaces on which the diagonal operators and characters
act are themselves direct integrals with respect to some fibered measure
of a measurable field of Hilbert spaces on $X$, the underlying Borel space
of the source.  

This last observation, together with the fact that the group of 2-arrows
under 1-dimensional composition is the semidirect product of the fiber
and base groups suggests that the representation theory of 2-groups may
have applications to classical group representation theory in addition
the the applications to topology and physics which motivated it. 

Observe in the case of $\frak P$ that, decomposing 
the representation into irreducible 
representations of the base group gives the same type of decomposition
considered in the case of Dirac's 
expansors for $SO(3,1)$ in \cite{CY.Lorentz}.  The connection
with Dirac's expansors is even more intimate:  an important family of
irreducible representations of $\frak P$ correspond to the orbits
of $SO(3,1)$ acting on $\widehat{M^4}$ (which may be identified
with Minkowski space), and are thus the energy levels for the system
of harmonic oscillators used by Dirac.  Thus the expansors are a 
quantization of the objects of our 2-category as conjectured in
\cite{CY.Lorentz}.

Finally, we must consider how to describe the 2-intertwiners.  A modification
assigns to each object of a bicategory a 2-arrow between the
1-arrows assigned by the source and target, such that the pillows

\begin{center}

\setlength{\unitlength}{4144sp}%
\begingroup\makeatletter\ifx\SetFigFont\undefined
% extract first six characters in \fmtname
\def\x#1#2#3#4#5#6#7\relax{\def\x{#1#2#3#4#5#6}}%
\expandafter\x\fmtname xxxxxx\relax \def\y{splain}%
\ifx\x\y   % LaTeX or SliTeX?
\gdef\SetFigFont#1#2#3{%
  \ifnum #1<17\tiny\else \ifnum #1<20\small\else
  \ifnum #1<24\normalsize\else \ifnum #1<29\large\else
  \ifnum #1<34\Large\else \ifnum #1<41\LARGE\else
     \huge\fi\fi\fi\fi\fi\fi
  \csname #3\endcsname}%
\else
\gdef\SetFigFont#1#2#3{\begingroup
  \count@#1\relax \ifnum 25<\count@\count@25\fi
  \def\x{\endgroup\@setsize\SetFigFont{#2pt}}%
  \expandafter\x
    \csname \romannumeral\the\count@ pt\expandafter\endcsname
    \csname @\romannumeral\the\count@ pt\endcsname
  \csname #3\endcsname}%
\fi
\fi\endgroup
\begin{picture}(2220,2641)(38,-1869)
\thinlines
\put(1351, 29){\oval(1530,362)[bl]}
\put(1351, 29){\oval(1530,362)[br]}
\put(2116, 29){\vector( 0, 1){0}}

\put(2116,-1141){\vector( 0,-1){0}}
\put(1351,-1141){\oval(1530,362)[tr]}
\put(1351,-1141){\oval(1530,362)[tl]}

\put(2123,382){\vector( 0,-1){0}}
\put(1358,382){\oval(1530,362)[tr]}
\put(1358,382){\oval(1530,362)[tl]}

\put(1343,-1486){\oval(1530,362)[bl]}
\put(1343,-1486){\oval(1530,362)[br]}
\put(2108,-1486){\vector( 0, 1){0}}

\put(2168, -1){\vector( 0,-1){990}}

\put(518,-24){\vector( 0,-1){990}}

\put(856,-331){\makebox(0,0)[lb]{\smash{\SetFigFont{12}{14.4}{it}$\Psi_*
$}}}
\put(1283,104){\makebox(0,0)[lb]{\smash{\SetFigFont{12}{14.4}{it}$\Downarrow
\phi_*
$}}}
\put(1306,-1306){\makebox(0,0)[lb]{\smash{\SetFigFont{12}{14.4}{it}$\Downarrow
\phi_*
$}}}
\put(931,-593){\makebox(0,0)[lb]{\smash{\SetFigFont{12}{14.4}{it}$\Phi_G
$}}}
\put(481,156){\makebox(0,0)[lb]{\smash{\SetFigFont{12}{14.4}{it}$X
$}}}
\put(2116,142){\makebox(0,0)[lb]{\smash{\SetFigFont{12}{14.4}{it}$Y
$}}}
\put( 38,-519){\makebox(0,0)[lb]{\smash{\SetFigFont{12}{14.4}{it}$R(G)
$}}}
\put(849,-1869){\makebox(0,0)[lb]{\smash{\SetFigFont{12}{14.4}{it}$\Psi_*
$}}}
\put(466,-1359){\makebox(0,0)[lb]{\smash{\SetFigFont{12}{14.4}{it}$X
$}}}
\put(2124,-1389){\makebox(0,0)[lb]{\smash{\SetFigFont{12}{14.4}{it}$Y
$}}}
\put(2258,-616){\makebox(0,0)[lb]{\smash{\SetFigFont{12}{14.4}{it}$S(G)
$}}}
\put(1343,-879){\makebox(0,0)[lb]{\smash{\SetFigFont{12}{14.4}{it}$\Psi_G
$}}}
\put(1808,637){\makebox(0,0)[lb]{\smash{\SetFigFont{12}{14.4}{it}$\Phi_*
$}}}
\put(1734,-893){\makebox(0,0)[lb]{\smash{\SetFigFont{12}{14.4}{it}$\Phi_*
$}}}
\put(1133,-421){\makebox(0,0)[lb]{\smash{\SetFigFont{12}{14.4}{it}$\Rightarrow
$}}}
\put(1411,-691){\makebox(0,0)[lb]{\smash{\SetFigFont{12}{14.4}{it}$\Rightarrow
$}}}
\end{picture}

\end{center}

\noindent commute.

Thus in our case, we need a single 2-arrow.  
Using the same sort of reasoning as in the case of representations and
1-intertwiners to restrict our attention to the behavior on partial
measurable line bundles, we can show the following:
If $\Phi$ and $\Psi$ are
1-intertwiners with $\Phi_* = \Phi_{{\cal K}, \{d\mu_y\}}$ and
$\Psi_* = \Phi_{{\cal L}, \{d\nu_y\}}$, a 2-intertwiner 
$\phi:\Phi \rightarrow \Psi$ is given by a bounded field of operators

\[ \phi_y:\int^\oplus_{x\in X} K_{<x,y>} 
 d\mu_y(x) \rightarrow
\int^\oplus_{x\in X} L_{<x,y>} 
d\nu_y(x) \]

\noindent such that the pillows commute.

Again, for any measure on $Y$ equivariant under the base group, taking
direct integrals gives a recognizable structure:  the 2-intertwiner
becomes an intertwiner in the ordinary sense between the representations
of the base group (only) arising by taking direct integrals of the
source and target 1-intertwiners.  This intertwiner is, however, given
by a diagonal operator with respect to the direct integral structure.

\section{Tensor Products}

As observed in \cite{Y.meas} the bicategory {\bf Meas} admits
a monoidal structure induced by the cartesian product of measure spaces.

\begin{defin}
Let $\odot:{\bf Meas}\times {\bf Meas}\rightarrow {\bf Meas}$ be given on
objects by $(X,S)\odot (Y,T) = (X\times Y, S*T)$, where $S*T$ denotes the
Borel structure induced by all products $A\times B$ for $A\in S$ and $B\in
T$.

$\odot$ is then defined on 1-arrows as follows:  for $F = \Phi_{{\cal K}
\{d\mu_\xi\}}$ and $G = \Phi_{{\cal L},\{d\nu_\theta\}}$
$F\odot G:{\bf Meas}(X)\odot 
{\bf Meas}(Y)\rightarrow {\bf Meas}(X^\prime)\odot {\bf Meas}(Y^\prime)$
is given on objects by

\[
{\cal H}\mapsto \int_{(x,y)\in X\times Y}^\oplus
{\cal H}_{(x,y)} \otimes {\cal K}_{<x,\xi>}
 \otimes {\cal L}_{<y,\theta>}  d\mu_\xi (x) \times
d\nu_\theta (y) \]

And by the same formula {\em mutatis mutandis} on arrows.

$\odot$ on 2-arrows $\phi:F\rightarrow F^\prime$ and $\psi:G\rightarrow 
G^\prime$, where $F = \Phi_{{\cal K}
\{d\mu_\xi\}}$, $F^\prime = \Phi_{{\cal K}^\prime 
\{d\mu^\prime_\xi\}}$,
$G = \Phi_{{\cal L},\{d\nu_\theta\}}$, and 
$G^\prime = \Phi_{{\cal L}^\prime,\{d\nu_\theta\}}$, is given
as follows:

From \cite{Y.meas} recall that a natural transformation $\phi$ between
measurable functors $F$ and $F^\prime$ as above is determine by a
measurable field of operators on $X^\prime$

\[ \{ \phi_\xi:\int_{x\in X}^\oplus {\cal K}d\mu_\xi(x)
\rightarrow \int_{x\in X}^\oplus {\cal K}^\prime d\mu^\prime_\xi(x)
\} \]

\noindent and similarly $\psi$ is given by a bounded field of
operators $\{ \psi_\theta \}$ on $Y^\prime$. $\phi \odot \psi$ is
then the natural transformation obtained by pre- and post-composing
$\{\phi_x \otimes \psi_y\}$ with the natural isomorphisms obtained
by applying the distributivity of $\otimes$ over direct integrals and
the categorical Fubini's Theorem of \cite{Y.meas}
\end{defin}

Just as the tensor product structure on {\bf VECT} induces a monoidal
structure on categories of representations of groups, so the monoidal
bicategory structure $\odot$ on {\bf Meas} induces a monoidal bicategory
structure on the bicategory of representations of any 2-group.

%The product $R\odot S$ of two representations 
%will always decompose into
%representations whose underlying Borel space is an orbit of the action
%of the base group on the product $R(*)\times S(*)$.  To more carefully
%formalize this, consider

%\begin{defin}
%A {\em $Y$-indexed measurable family of representations of an
%automorphic 2-group $\frak G$} is
%a representation of $\frak G$ with underlying Borel space $X\times Y$ such
%that the action of the base-group ${\frak G}_0$ stablizes the sets 
%$\{x\}\times Y$ for all $x\in X$.
%\end{defin}

\section{Additive Reflections as Projections}

It appears that the natural analogue at the level of 1-arrows for the 
projections occurring as labels of faces in the state-sum constructions
in dimension 3 are additive reflections, in particular those induced
by bimeasurable inclusions.

Recall that a reflection in the categorical sense is a pair of adjoint
functors $i \vdash r$ such that the unit $\eta: r(i)\rightarrow 
Id_{{\rm source}(i)} = Id_{{\rm target}(r)}$ is an isomorphism.  We
call a reflection between additive categories, both of whose functors
are additive functors {\em an additive reflection}.

Recall from \cite{Y.meas} that a function between measurable spaces
is bimeasurable if both direct and inverse image preserve measurable sets.
A bimeasurable inclusion $i:(X,S)\rightarrow (Y,T)$
induces an additive reflection $i_*:{\bf Meas}(Y,T)\rightarrow 
{\bf Meas}(X,S) \vdash  i^*{\bf Meas}(X,S)\rightarrow {\bf Meas}(Y,T)$.

Additive reflections of this form provide a kind of decomposition of
tensor products of representations of any automorphic 2-group:
Consider two representations $R$ and $R^\prime$ with underlying
Borel spaces $(X,S)$ and $(X^\prime, S^\prime)$.  

Now, consider the representation $R \odot R^\prime$.  Its underlying
Borel space is $(X\times X^\prime, S*S^\prime)$, and the base group
${\frak G}_0$ acts by the diagnonal action, and the fiber group acts
by the (tensor) product of the characters assigned to the coordinates
on the fibers of measurable 
fields.

The action of ${\frak G}_0$ decomposes $X\times X^\prime$ into orbits,
each of which inherits a Borel structure from the product, and thus
is included by a bimeasurable inclusion.

Now it is clear that representations in which the base group acts
transitively on the underlying Borel space are irreducible in the sense
of admitting no proper non-empty subobjects. (Note, the linear structure
here exists at the level of 1- and 2-arrows, not objects, so the
minimal subrepresention is not a ``0-representation'' but the one
with empty underlying Borel space.)

Additive relections
along inclusions of orbits
thus provide a ``decomposition'' of any representation into irreducibles.
Moreover, both functors in these additive reflections are measurable:
each is represented by a measurable line bundle on the product concentrated
on the graph of the inclusion and direct integration with respect to
measures concentrated on the image (or preimage) of the point.

\section{Coloring triangulations} \label{ss}

Although it is beyond the scope of the present work to give complete 
constructions of either 4-dimensional TQFT's or the model for quantum 
gravity which
are the principal motivations for considering representations of 2-groups
in {\bf Meas}, we wish to give an indication of the conjectural constructions
for which we have developped this theory.  

Both proceed by ``coloring'' the simplexes of 
a triangulated 4-manifold (considered as a space-time in the QG case)
with objects, 1-arrows, and 2-arrows of ${\bf Rep}({\frak G})$ (for
${\frak G} = {\frak P}$ in the QG case), organizing the colorings into
a suitable measure space, and taking the volume of the measure space.  

Notice that the more familiar state-sum constructions for TQFT's, and
the construction for Euclidean quantum gravity proposed by Barrett and Crane
\cite{BC1} can be phrased in these terms:  the product of 6j- or 10j-symbols
and (quantum) dimensions (and their reciprocals) is defining a discrete
measure on the space of colorings.  

In the present case, the space of colorings will not admit a discrete measure
with the appropriate properties.

At least two versions of the construction readily come to mind.  Both 
at least formally lead to TQFT's when applied to any automorphic 2-group,
and either could, in the case of the Poincar\'{e} 2-group ${\frak P}$,
lead to a good quantization of general relativity provided a suitable
measure, concentrated on colorings which embody a quantum analog of the
geometric restriction that a face as a bivector is the wedge of any
two of its edges as vectors, can be discovered.

In the first, edges are colored with those representations which have the base
group with its natural Borel structure as underlying Borel space. These
are indexed by equivariant families of characters, which indexing set 
inherits a Borel structure from the Borel structure used to define
Plancherel measure on the dual of the (abelian) fiber group.

Face labels are then the additive reflections from one edge label to
the tensor product of the other two, and are thus indexed by a family
of orbits in the product ${\frak G}_0 \times {\frak G}_0$ with the diagonal
action, namely those orbits for which product of a character
in the equivariant family giving one tensorand with a character in the
family giving the other tensorand lies in the family giving the
label on the remaining edge.  As the family of all 
the orbits of the diagonal action is indexed by the base group,
and thus has a natural Borel structure, the family indexing the
admissible face labels inherits a Borel structure.      

In the second, edge are colored with those representation which have
as underlying Borel space an orbit of the dual group $\hat{E}$ of the fiber
group $E$ under the action of the base group.  These orbits inherit a
Borel structure from that used in defining the Plancherel measure on 
$\hat{E}$.  Two variations are possible, one in which all such
representations are permitted, and another in which we require
that only the ``tautologous'' representation, in which
at each point the Hilbert space is 1-dimensional and has the fiber group
acting by the character which names the point, be used.  
For physical applications, this latter seems most promising as the
permitted colorings are precisely the mass shells.

\begin{flushleft}
\bibliography{MeasCatand2Groups}
\end{flushleft}

\end{document}